\documentclass[a4paper, 11pt]{article}
\usepackage[english]{babel}
\usepackage{latexsym}
\usepackage{amsfonts, amsthm}
\usepackage{amssymb}
\usepackage{amscd}
\usepackage{amsmath}
\usepackage{graphics}
\usepackage{eepic}
\usepackage{epsfig}
\usepackage[noxcolor]{pstricks}
\newtheorem{theorem}{Theorem}
\newtheorem{lemma}[theorem]{Lemma}
\newtheorem{proposition}[theorem]{Proposition}
\newtheorem{corollary}[theorem]{Corollary}

\newtheorem{conjecture}[theorem]{Conjecture\/}

\newcommand{\C}{{\mathbb C}}

\newcommand{\fie}{\varphi}
\newcommand{\I}{\mathcal{I}}
\newcommand{\<}{\backslash}

\begin{document}

\title{\textbf{The holomorphy conjecture for ideals in dimension two}}
\author{Ann Lemahieu and Lise Van Proeyen \footnote{Ann Lemahieu, Lise Van Proeyen; K.U.Leuven, Departement Wiskunde,
Celestij\-nenlaan 200B, B-3001 Leuven, Belgium, email:
lemahieu@mathematik.uni-kl.de, lise.vanproeyen@wis.kuleuven.be. The
research was partially supported by the Fund of Scientific Research
- Flanders (G.0318.06) and MEC PN I+D+I MTM2007-64704.}
\date{}}
\maketitle {\footnotesize \emph{\noindent \textbf{Abstract.---} The
holomorphy conjecture predicts that the local Igusa zeta function
associated to a hypersurface and a character is holomorphic on
$\mathbb{C}$ whenever the order of the character does not divide the
order of any eigenvalue of the local monodromy of the hypersurface.
In this note we propose the holomorphy conjecture for arbitrary
subschemes at the level of the topological zeta function and we
prove this conjecture for subschemes defined by an ideal that is
generated by a finite number of complex polynomials in two
variables.}}
\\ \\
 ${}$ \begin{center}
\textsc{0. Introduction}
\end{center} ${}$\\
During the last decades a lot of research has been done on the poles
of zeta functions. These zeta functions are associated to a
polynomial or to a germ of a holomorphic function. The zeta
functions we study, for instance the Igusa zeta function, the
motivic zeta function and the topological zeta function, are
rational functions that can be computed from an embedded resolution
of singularities. We will work locally, say around the origin. Every
irreducible component of the total transform (through the embedded
resolution) of the germ at the origin gives rise to a candidate pole
of the zeta function. It is striking that a lot of candidate poles
are cancelled. For plane curves one can easily determine the poles
from the resolution graph: in \cite{polen dim2} Veys shows that a
candidate pole is a pole if and only if it is associated to an
exceptional component intersecting at least three times other
components or to an irreducible component of the strict transform.
Notice also the similarity with jumping numbers of multiplier ideals
for curves (see
\cite{Smith} and \cite{kevin}) and stable curves. \\
\indent These poles also play a main role in some very mysterious
conjectures. For instance the monodromy conjecture predicts that a
pole of a local zeta function gives rise to an eigenvalue of the
local monodromy. This conjecture has been proven for several cases
(see for example \cite{Igusa}, \cite{Loeser num data},
\cite{loeser2}, \cite{luengo}, \cite{luengo2},
\cite{monodconjtoric}, \cite{polen dim2} and \cite{VeysPVI}) but is
still open in general and until now one does not understand
why these poles and eigenvalues would be related in such a way. \\

\indent The holomorphy conjecture, proposed by Denef in \cite{Denef
Bourbaki rapport}, essentially states that the Igusa zeta function
associated to a polynomial and a character is holomorphic on
$\mathbb{C}$ when the order of the character is no divisor of the
order of any eigenvalue of monodromy. There exists an analogous
version of the holomorphy conjecture for the topological zeta
function associated to a polynomial and an integer (the integer
corresponds to the order of the character).

In this article we introduce the holomorphy conjecture for ideals in
$\C[x_1,\ldots,x_n]$. In particular we show that the topological
zeta function associated to an arbitrary subscheme of $\C^n$ and an
integer is well-defined. The notion of embedded resolution is here
replaced by the notion of log-principalisation of the ideal defining
the subscheme. In Section $2$ we go on by providing some preliminary
results in dimension $2$ which we will use in Section $3$ to prove
the holomorphy conjectures for ideals in $\mathbb{C}[x,y]$.
\\ \\
 ${}$ \begin{center}
\textsc{1. The holomorphy conjecture}
\end{center}
${}$\\
The holomorphy conjecture has to do with local zeta functions being
holomorphic on the whole complex plane and the orders of the
eigenvalues of local monodromy. Verdier introduced a notion of
eigenvalues of monodromy for ideals, coinciding with the classical
notion for principal ideals (see \cite{verdier}). Based on this
notion of Verdier, the second author and Veys gave a criterion \`{a}
la A'Campo for being an eigenvalue of monodromy of a given ideal.

To recall this criterion, fix an ideal $\I \subset
\C[x_1,\ldots,x_n]$. Let $Y$ be the zero locus of $\I$ in $X:=\C^n$,
containing the origin $0$. We denote the blowing-up of $X$ with
respect to $\I$ by $Bl_\I X$. Now consider a log-principalisation
$\psi: \tilde{X} \rightarrow X$ of $\I$ (the existence of that is
guaranteed by Hironaka in \cite{hironaka}). This means that $\psi$
is a proper birational map from a nonsingular variety $\tilde{X}$
such that the total transform $\psi^*\I$ is locally principal and
moreover is the ideal of a simple normal crossings divisor. Let
$\sum_{i \in S} N_i E_i$ denote this divisor, written in such a way
that the $E_i, i \in S$, are the irreducible components occurring
with multiplicity $N_i$. Let $\nu_i-1$ be the multiplicity of $E_i$
in the divisor $\psi^*(dx_1 \wedge \cdots \wedge dx_n)$. The couples
$(N_i,\nu_i), i \in S$, are called the numerical data of the
log-principalisation $\psi$. For $I \subset S$, denote $E_I:=\cap_{i
\in I}E_i$ and $E_I^\circ :=E_I \< (\cup_{j \in S, j \notin I}E_j)$.
We denote furthermore the topological Euler-Poincar\'e
characteristic by $\chi(\cdot)$. By the Universal Property of
Blowing Up, there exists a unique morphism $\fie$ that makes the
following diagram commutative.
\begin{center}
\begin{pspicture}(0,0)(3.2,3.2)
\rput(0,3){$\tilde{X}$} \rput(3,3){$Bl_\I X$} \rput(3,0){$X$}
\psline{->}(0.3,2.7)(2.7,0.3) \psline{->}(3,2.7)(3,0.3)
\psline{->}(0.3,3)(2.5,3) \rput(1.5,1){$\psi$}
\rput(3.2,1.5){$\pi$}\rput(1.5,3.2){$\fie$}
\end{pspicture}
\end{center}

\noindent The eigenvalues of monodromy defined by Verdier were
originally associated to points of the normal cone $C_Y X$. These
eigenvalues are the same for points on the same ruler. Hence, we can
define the eigenvalues of monodromy for every point of the
projectivisation of the normal cone, which is the support $E$ of the
total transform of $\I$ by $\pi$ in $Bl_\I X$. We will use the
following criterion by the second author and Veys (see \cite{VV
monodromie}).

\begin{theorem} \label{criterion}
The number $\alpha$ is an eigenvalue of monodromy of $\I$ if and
only if there exists a point $e \in E$ such that $\alpha$ is a zero
or pole of the function
\[Z_{\I,e} (t) =
\prod_{j \in S} (1-t^{N_j})^{\chi(E_j^\circ \cap \varphi^{-1}(e))}.
\]
\end{theorem}

To a positive integer $d$ and the ideal $\I$ one can associate a
local topological zeta function at the point $0$ which is a rational
function in one complex variable.
\[Z_{top,\I}^{(d)}(s):=\sum_{\substack{I \subset S \\ d \mid N_i}} \chi (E_I^{\circ}
\cap \psi^{-1}\{ 0 \}) \prod_{i \in I} \frac{1}{N_i s+ \nu_i}.\]

\noindent When the ideal $\I$ is principal, then $Z_{top,\I}^{(d)}$
can be defined as a limit of $p$-adic Igusa zeta functions with
character of order $d$ and hence is well-defined (see
\cite{DenefLoeser1}). We now show that these topological zeta
functions are also well-defined (independent of the chosen
log-principalisation) for ideals.

Given two log-principalisations $\psi_1: \tilde{X_1} \rightarrow
\C^n$ and $\psi_2: \tilde{X_2} \rightarrow \C^n$ for $\I$, then
there exists a birational map from $\tilde{X_1}$ to $\tilde{X_2}$
making the diagram commutative. The Weak Factorization Theorem (see
\cite{AKMW-factorization}) states that this birational map can be
factorized as a series of consecutive blowing-ups and blowing-downs
with smooth centra that have normal crossings with the total
transform at the intermediate stage. To assure that
$Z_{top,\I}^{(d)}$ is well-defined, the Weak Factorization Theorem
also guarantees that it is sufficient to show that a
log-principalisation $\psi_1 \circ \sigma$ yields the same function
$Z_{top,\I}^{(d)}$ as $\psi_1$ does, with $\sigma: \tilde{X_1'} \to
\tilde{X_1}$ a blowing-up with smooth centre $C$ having normal
crossings with $\psi_1^*\I$.

For a fixed point $p \in C$, we check that the contribution of $p$
to the topological zeta function coming from $\psi_1$ is equal to
the contribution of $\sigma^{-1}(p)$ to the topological zeta
function coming from  $\psi_1 \circ \sigma$. We will not give all
details here, but we provide the main ideas. Let $E_s$ be the
exceptional divisor of $\sigma$. Then $E_s$ is a projective bundle
over $C$ and $\sigma^{-1}(p) \cong \mathbb{P}^{\textnormal{
codim}(C)-1}$. Suppose $\{1, \ldots, k\}$ are the indices of the
components of $\psi_1^*\I$ that contain $C$ and $\{k+1, \ldots ,
l\}$ the indices of the components that intersect $C$ transversally
and which contain $p$.
For $i\in \{1,\ldots,l\}$, we put $E_i(N_i,\nu_i)$ for the
corresponding components  in $(\psi_1 \circ \sigma)^*\I$ with their
numerical data. As $N_s = \sum_{i=1}^k N_i$, notice that it is
impossible that at the same time $d \mid N_s$ and that there exists
exactly one of the $N_i (1 \leq i \leq k)$ that is no multiple of
$d$.

Then we need the Euler-Poincar\'e characteristics of the $E_I^\circ
\cap \sigma^{-1}(p)$ for sets $I$ that contain $s$ and a subset of
$\{1, \ldots , l\}.$ These are always zero if $\{k+1, \ldots l\}
\not \subset I.$ If $\{s, k+1, \ldots l\} \subset I$, we define
$m=\#(I\cap \{1, \ldots , k\})$. Then
\begin{eqnarray*}
\chi(E_I^\circ \cap \sigma^{-1}(p)) & = & 0 \qquad \qquad\qquad\
\mbox{ if } m
\leq k-2, \\
& = & 1 \qquad\qquad \qquad \ \mbox{ if } m = k-1, \\& = &
\mbox{codim}(C) - k \quad \mbox{ if } m = k.
\end{eqnarray*}
We can use this to show that the mentioned contributions are zero if
$d$ is no divisor of some $N_i, i=1,\ldots,l.$ Otherwise the
contributions are equal to
$\left(\prod_{i=1}^l(N_is+\nu_i)\right)^{-1}.$
\\ \\
In this article we will prove the holomorphy conjecture in the
special case that $\I$ is an ideal in $\C[x,y]$.

\begin{conjecture} (Holomorphy Conjecture)\\
Let $d$ be a positive integer. If $d$ does not divide the order of
any eigenvalue of monodromy associated to the ideal $\I$ in points
of $\pi^{-1}\{0\}$, then $Z_{top,\I}^{(d)}$ is holomorphic on the
complex plane.
\end{conjecture}

\noindent In the statement of this conjecture we use the word
`holomorphic', but actually we are going to show that the mentioned
function is identically zero. This terminology has its origins in
the context of $p$-adic Igusa zeta functions.

When the ideal is principal, then this conjecture has been shown by
Veys in \cite{Veys holomorphy} for plane curves. Veys and the first
author confirmed the conjecture for surfaces that are general for a
toric idealistic cluster (see \cite{monodconjtoric}). We will show
the following statement: \\ \\ \emph{Let $d$ be a positive integer
that does not divide the order of any eigenvalue of monodromy in
$\pi^{-1}\{0\}$ associated to the ideal $\I \subset
\mathbb{C}[x,y]$. The exceptional components $E_i$ for which $d$
divides $N_i$ satisfy:
\begin{enumerate}
\item $\chi(E_i^{\circ})=0$;
\item these components do not intersect.
\end{enumerate}}
\noindent We use this result to prove the holomorphy conjecture for
ideals in $\mathbb{C}[x,y]$. The structure of our proof is inspired
by the structure of the proof by Veys in \cite{Veys holomorphy} for
plane curves. This proof was actually given for the local $p$-adic
Igusa zeta function with character associated to a hypersurface.
However, in the context of ideals a definition of Igusa zeta
function with character is not (yet) known. \vspace{\baselineskip}

\noindent \textbf{Example.} We consider the ideal $\I=(x^2y^4,
x^{34}, y^6) \subset \C[x,y].$ A log-principalisation of $\I$
consists of eight successive blowing-ups. The intersection diagram
together with the associated numerical data can be found in the
following figure.
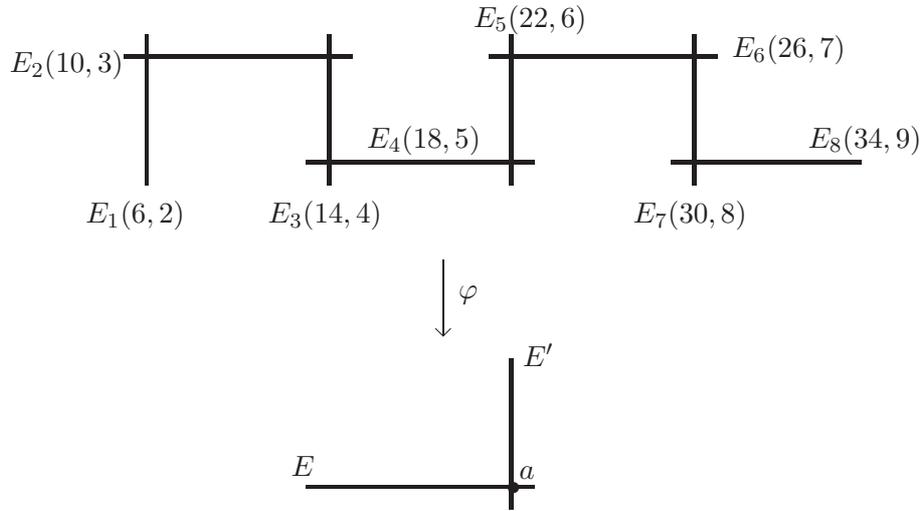
\begin{figure}[h]
\unitlength=1mm
\begin{center}
\begin{picture}(90,65)(0,-40)
\linethickness{0.50mm}

\put(3,0){\line(0,1){20}} \put(0,17){\line(1,0){30}}
\put(27,0){\line(0,1){20}} \put(24,3){\line(1,0){30}}
\put(51,0){\line(0,1){20}} \put(48,17){\line(1,0){30}}
\put(75,0){\line(0,1){20}} \put(72,3){\line(1,0){25}}

\put(-5,-5){$E_1(6,2)$} \put(-15,15){$E_2(10,3)$}
\put(19,-5){$E_3(14,4)$} \put(32,5){$E_4(18,5)$}
\put(46,21){$E_5(22,6)$} \put(80,17){$E_6(26,7)$}
\put(67,-5){$E_7(30,8)$} \put(90,5){$E_8(34,9)$}

\linethickness{0.10mm} \put(42,-10){\line(0,-1){10}}
\put(42,-20){\line(-1,1){1}}\put(42,-20){\line(1,1){1}}

 \put(44,-15){$\fie$}

 \linethickness{0.50mm} \put(24,-40){\line(1,0){30}}
\put(51,-43){\line(0,1){20}} \put(22,-38.5){$E$}
\put(52.5,-24){$E'$}

\put(50.3,-41){$\bullet$} \put(52,-39){$a$}

\end{picture}
\end{center}
\end{figure}

\noindent We use Theorem \ref{criterion} to find the eigenvalues of
monodromy. The exceptional curves $E_2, \ldots, E_7$ are contracted
by the map $\fie$ to the intersection point $a$ of the exceptional
components $E$ and $E'$ in $Bl_\I \,\C^2.$ The exceptional curves
$E_1$ and $E_8$ are respectively mapped surjectively to $E$ and
$E'$. As eigenvalues of monodromy we get the 6th roots of unity and
the 34th roots of unity. For instance $d=5$ is no divisor of the
order of an eigenvalue of monodromy. The components $E_2$ and $E_7$
satisfy $\chi(E_2^\circ)= \chi(E_7^\circ)=0$ and have an empty
intersection. This implies that $Z_{top, \I}^{(5)}(s)$ is equal to
zero. \qed
\\ ${}$ \begin{center}
\textsc{2. Preliminary results}
\end{center} ${}$\\
\noindent From now on, we consider an ideal $\I:=(f_1, \ldots, f_r)$
in $\C[x,y]$. Notice that a log-principalisa\-tion of an ideal also
gives an (non-minimal) embedded resolution for all members of some
Zariski open subset of the linear system $\{\lambda_1 f_1 + \cdots +
\lambda_r f_r = 0 \, | \, \lambda_1, \ldots, \lambda_r \in \C \}$.
We will call the elements for which the principalisation gives an
embedded resolution \emph{totally general for $\I$}. Moreover, the
numerical data associated to the principalisation and to the
embedded resolution are the same. A proof of this statement can be
found in \cite[\S 2]{polen ideaal dim2}. Let us write $\I$ as $\I =
(h)(f_1', \ldots , f_r')$ with $(f_1', \ldots , f_r')$ finitely
supported. We will say that a totally general element for $(f_1',
\ldots , f_r')$ is \emph{general for $\I$}.
We will use the notation introduced in Section $1$. In particular
the $E_i, i \in S$, will be the irreducible components of
$\psi^{-1}\I$. We choose a totally general element $f$ for $\I$ and
we can write $\psi^{-1}(f^{-1}\{0\}) = \sum_{i \in T}N_i E_i$, with
$T$ a set containing $S$. Let $k_i, i \in S,$ be the number of
intersection points of $E_i$ with other components of $\psi^{-1}\I$.
Analogously, for $i \in T,$ let $k_i'$ be the number of intersection
points of $E_i$ with other components of $\psi^{-1}(f^{-1}\{0\}).$
So $k_i \leq k_i'$ for $i \in S,$ with equality if and only if $E_i$
is not intersected by the strict transform of a general element for
$\I$. We will use the following congruence.
\begin{lemma}\emph{\cite[Lemme II.2]{Loeser num data}} \label{lemma Loeser relatie N_i} If we fix one
exceptional curve $E_i,$ intersecting $k_i'$ times other components
$E_1, \ldots, E_{k_i'}$ of the embedded resolution
$\psi^{-1}(f^{-1}\{0\})$, then
$$\sum_{j=1}^{k_i'} N_j \equiv 0 \mbox{ mod }N_i.$$
\end{lemma}
\noindent Veys shows the following result in his proof for the
holomorphy conjecture for plane curves. He proved this for the
minimal embedded resolution, but the proof remains valid for
non-minimal resolutions induced by log-principalisations.
\begin{lemma} \emph{\cite[Lemma 2.3]{Veys holomorphy}} \label{lemma veys E_0 snijdt 1 keer}Let $E_0$ be an
exceptional curve with $k_0' = 1.$ Then for some $r \geq 1$ there
exists a unique path
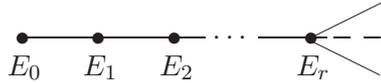
\begin{figure}[h]
\unitlength=1mm
\begin{center}
\begin{picture}(40,5)(0,3)
\linethickness{0.15mm}

\put(1,5){\line(1,0){24}} \put(26,4){$\cdots$}
\put(32,5){\line(1,0){6}}

\put(39,5){\line(2,1){10}} \put(38,5){\line(2,-1){10}}
\put(38,5){\line(1,0){3}} \put(42,5){\line(1,0){2}}
\put(46,5){\line(1,0){2}}

\put(0,3.9){$\bullet$} \put(10,3.9){$\bullet$}
\put(20,3.9){$\bullet$} \put(38,3.9){$\bullet$}

\put(-1,0){$E_0$} \put(9,0){$E_1$} \put(19,0){$E_2$}
\put(37,0){$E_r$}

\end{picture}
\end{center}
\end{figure}

\noindent in the resolution graph consisting entirely of
exceptional curves, such that \begin{enumerate}
\setlength{\itemsep}{-3pt} \item $k_j' = k_j =2$ for $j=1, \ldots,
r-1;$ \item $k_r' \geq 3;$ \item $N_0 | N_j$ for all $j=1, \ldots
, r;$ \item $N_0 < N_1 < \cdots < N_r.$
\end{enumerate}
\end{lemma}
\noindent We will now provide a set of eigenvalues of monodromy.
Recall that the Rees components of an ideal $\I$ are the irreducible
components of the exceptional divisor on $\overline{Bl_\I X}$. Let
$n:\overline{Bl_\I X} \rightarrow Bl_\I X$ be the normalization map
and let $\sigma: \tilde{X} \rightarrow \overline{Bl_\I X}$ be such
that $\fie=n \circ \sigma$. We will also call the corresponding
exceptional components in $\tilde{X}$ Rees components, so an
exceptional component $E$ in $\tilde{X}$ is Rees if and only if
dim$(\sigma(E))=$dim$(E)$. As the normalization map is a finite map,
being contracted by $\fie$ is equivalent to being contracted by
$\sigma$. Theorem \ref{criterion} gives us:
\begin{corollary} \label{lemma rees dan eigenwaarde}
If the exceptional component $E$ in $\tilde{X}$ is Rees for $\I$,
then all $N$th roots of unity are eigenvalues of monodromy.
\end{corollary}
\noindent We can recognize these Rees components in the resolution
graph in a very easy way.
\begin{lemma}\label{lemmaeigenwaarde}
An irreducible component $E$ on $\tilde{X}$ is contracted by the map
$\fie: \tilde{X} \rightarrow Bl_\I \,\mathbb{C}^2$ if and only if
the strict transform of a general element for $\I$ does not
intersect $E$.
\end{lemma}

\noindent \emph{Proof.} Let $D$ be the Cartier divisor on
$\tilde{X}$ such that
$\I\mathcal{O}_{\tilde{X}}=\mathcal{O}_{\tilde{X}}(-D)$ and let $F$
be the Cartier divisor on $Bl_\I \,\C^2$ such that
$\I\mathcal{O}_{Bl_\I \,\C^2}=\mathcal{O}_{Bl_\I \,\C^2}(-F)$. Then
by the projection formula one has $(-D)\cdot E = -\fie^*(F)\cdot E =
(-F)\cdot \fie_*E$.

Suppose $E$ is contracted by $\fie$, then $\fie_*E=0$ and $(-D)\cdot
E=0$. If $E$ is not contracted by $\fie$, then $\fie_*E=k\fie(E)$
for some strictly positive integer $k$. If $\I=(f_1,\ldots,f_r)$,
then we have a map $j: Bl_\I \,\C^2 \stackrel{i}{\hookrightarrow}
\C^2 \times \mathbb{P}^{r-1} \stackrel{pr}{\rightarrow}
\mathbb{P}^{r-1}$ where $i$ is the canonical embedding of $Bl_\I
\,\C^2$ in $\C^2 \times \mathbb{P}^{r-1}$ and $pr$ is the projection
map. We have $\I\mathcal{O}_{Bl_\I \,\C^2}=j^* \mathcal{O}(1)$ or
$-F=j^*H$ in Pic $Y$ for $H$ a divisor class of hyperplanes on
$\mathbb{P}^{r-1}$.

Then again by the projection formula we get $(-F)\cdot \fie(E)=H
\cdot j_*\fie(E)=\textnormal{deg}(\fie(E)) > 0$ and thus $(-D)\cdot
E
> 0$.


We now write $\I=h\I'$ with $\I'$ an ideal of finite support. For a
totally general element $f=hf'$ for $\I$, we can write its total
transform $\psi^{-1}(f^{-1}\{0\})=D+S$, where $S$ is the strict
transform of $f'$. By the projection formula, one always has that
$(D+S)\cdot E=0$.
\\ \indent Combining these formulas, one gets the statement of Lemma \ref{lemmaeigenwaarde}.

\qed
\begin{proposition} \label{prop k_i>=3 dan eigenwaarde}
Let $E_j$ be an exceptional curve with $k_j \geq 3.$ Then $N_j$
divides the order of an eigenvalue of monodromy of $\I.$
\end{proposition}

\noindent \emph{Proof.} If $E_j$ is Rees for $\I$, then Corollary
\ref{lemma rees dan eigenwaarde} yields exactly this result. Suppose
now that $E_j$ is not Rees for $\I$ and let $a$ be the point on the
exceptional locus of $Bl_{\I}\,\C^2$ such that $a = \fie(E_j),$
where $\fie : \tilde{X} \to Bl_{\I}\,\C^2$. We define $S_a$ as the
set of indices $i \in S$ which satisfy $\fie(E_i) = a$. We see that
it is enough to prove that
$$\sum_{i \in S_a , N_j | N_i} \chi(E_i^\circ) \neq 0. $$
It is given that $\chi(E_j^\circ) < 0.$ We will now prove that every
positive contribution to this sum is cancelled by another negative
contribution. \\ \indent Suppose $\ell \in S_a, N_j | N_\ell$ and
$\chi(E_\ell^\circ) > 0$. This means that $\chi(E_\ell^\circ) =1$
and $k_\ell$ is equal to $1.$ If $k_\ell' \neq 1,$ then $E_\ell$ is
Rees for $\I$ and $N_j$ is a divisor of the order of an eigenvalue
of monodromy (Corollary \ref{lemma rees dan eigenwaarde}). If
$k_\ell' = 1$, then we can use Lemma \ref{lemma veys E_0 snijdt 1
keer} to find $E_r$ with $k_r' \geq 3.$

If $E_r$ is Rees for $\I$, Corollary \ref{lemma rees dan
eigenwaarde} tells us that $e^{\frac{2\pi i}{N_r}}$ is an eigenvalue
of monodromy and as $N_j|N_r,$ also $N_j$ divides the order of it.
Suppose now that $E_r$ is not Rees for $\I$. Then $E_{\ell+1},
\cdots,E_r$ are all contracted to the point $a$. Moreover we then
have $k_r=k_r'$ and thus $\chi(E_r^\circ) < 0.$  Now $N_j | N_\ell$
and $N_\ell | N_r$, so we have found a negative contribution
cancelling $\chi(E_\ell^\circ).$ \\ \indent We now check whether
there can exist two exceptional curves $E_\ell$ and $E_{\ell'}$ with
$\fie(E_\ell) = \fie(E_{\ell'}) = a, \chi(E_\ell^\circ) =
\chi(E_{\ell'}^\circ) = 1, N_j | N_\ell$ and $N_j | N_{\ell'},$ for
which the respectively associated $E_r$ and $E_{r'}$ are equal. We
know that $E_r$ is created later in the principalisation process
than $E_\ell, \ldots , E_{r-1}, E_{\ell'}, \ldots, E_{r'-1}.$ So at
the stage where $E_r$ is created, the resolution graph looks as
follows.
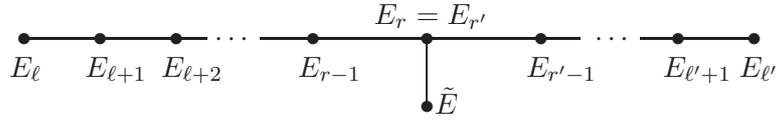
\begin{center}
\begin{picture}(100,15)(0,-4)
\linethickness{0.15mm} \put(1,5){\line(1,0){24}}
\put(26,4){$\cdots$} \put(32,5){\line(1,0){42}} \put(76,4){$\cdots$}
\put(82,5){\line(1,0){15}}

\put(0,4){$\bullet$} \put(10,4){$\bullet$} \put(20,4){$\bullet$}
\put(38,4){$\bullet$} \put(53,4){$\bullet$} \put(68,4){$\bullet$}
\put(86,4){$\bullet$} \put(96,4){$\bullet$}

\put(-1,0){$E_\ell$} \put(9,0){$E_{\ell+1}$}
\put(19,0){$E_{\ell+2}$} \put(37,0){$E_{r-1}$}
\put(47,7){$E_r=E_{r'}$} \put(67,0){$E_{r'-1}$}
\put(85,0){$E_{\ell'+1}$} \put(95,0){$E_{\ell'}$}

\put(53.9,5){\line(0,-1){9}}

\put(53,-5){$\bullet$} \put(55,-5){$\tilde{E}$}

\end{picture}
\end{center}

\noindent Note that by the principalisation process it is impossible
to have more than two exceptional curves intersecting $E_r.$ We
denote by $\tilde{E}$ the components of the strict transform of the
curves that belong to the support of $\I.$ These components might be
singular and are only present in the principalisation graph if $\I$
is not finitely supported. Since the principalisation graph is
connected, there are no other components at that moment. Hence, as
$E_j$ intersects at least three times, it follows that $N_j \geq
N_r$ and thus $N_j > N_\ell.$ This contradicts $N_j | N_\ell.$ \qed
\\ \\
 ${}$ \begin{center}
\textsc{3. Holomorphy conjecture for ideals in $\C[x,y]$}
\end{center} ${}$\\

\noindent Now we prove the holomorphy conjecture for the local
topological zeta function associated to an ideal in dimension two.

\begin{theorem}
Let $\I$ be an ideal in $\C[x,y]$ and $\pi: Bl_\I \,\C^2 \to \C^2$
be the blowing-up of $\C^2$ in the ideal $\I$. Suppose $d$ is a
positive integer that does not divide the order of any eigenvalue of
monodromy associated to the ideal $\I$ in points of $\pi^{-1}\{0\}$.
Then $Z_{top,\I}^{(d)}$ is identically $0$ on the complex plane.
\end{theorem}

\noindent \emph{Proof.} We search for components that contribute to
the local topological zeta function. If $\I$ is a principal ideal,
then we refer to \cite{Veys holomorphy}.

Suppose that $E_i(N_i, \nu_i)$ is an exceptional component of the
principalisation satisfying $d | N_i.$ By Corollary \ref{lemma rees
dan eigenwaarde} it follows that $E_i$ is not Rees for $\I$ and thus
$k_i=k_i'.$ If $k_i \geq 3,$ we use Proposition \ref{prop k_i>=3 dan
eigenwaarde} to see that $d$ would be a divisor of the order of a
monodromy eigenvalue.
If $k_i'=1,$ we use Lemma \ref{lemma veys E_0 snijdt 1 keer} to find
an exceptional curve $E_r$ with $k_r' \geq 3.$ If $k_r = k_r',$ we
are again in the situation of Proposition \ref{prop k_i>=3 dan
eigenwaarde}. Since $d |N_i$ and $N_i | N_r,$ this leads to a
contradiction. If $k_r \neq k_r',$ the component $E_r$ is Rees for
$\I$ and Corollary \ref{lemma rees dan eigenwaarde} brings the same
conclusion. Hence, we obtain that having $d|N_i$ for an exceptional
component $E_i(N_i, \nu_i)$ implies that $k_i = 2.$

Suppose now that $E_i(N_i, \nu_i)$ is a component of the support of
the weak transform satisfying $d|N_i.$ The only possible
contribution of $E_i$ comes from an intersection point of $E_i$ with
an exceptional component $E_j(N_j, \nu_j)$ for which $d|N_j.$ By
Corollary \ref{lemma rees dan eigenwaarde} it follows that $E_j$ is
not Rees for $\I$. Then we showed that there exists exactly one
other component $E_k$ that intersects $E_j$. From Lemma \ref{lemma
Loeser relatie N_i} it follows that $d|N_k$. If $E_k$ is Rees for
$\I$, then we have a contradiction. If $E_k$ is a component of the
support of the weak transform, then there is no Rees component in
the principalisation graph. This implies that $\I$ is a principal
ideal. If $E_k$ is exceptional and not Rees for $\I$, we can iterate
this argument. By finiteness of the resolution graph we should once
meet a component that is Rees for $\I$ or that is a component of the
support of the weak transform. This has been discussed before.

The only contribution to the topological zeta function can come from
an exceptional component $E_i$ with $\chi(E_i^\circ)=0$. In
particular, the contribution has to come from intersections with
other exceptional components. Suppose that $E_j$ is a component that
intersects $E_i$ and that $d|N_j.$ Then $E_j$ must be exceptional.
We do the same reasoning for $E_j$ and we find that $k_j$ must be
two. Suppose $E_k$ is the other component that intersects $E_j.$ By
Lemma \ref{lemma Loeser relatie N_i} we know that $d$ must divide
$N_k.$ We iterate this argument and get the existence of a component
$E(N,\nu)$ that is Rees for $\I$ and for which $d|N$. This
contradicts the choice of $d$ (Corollary \ref{lemma rees dan
eigenwaarde}) and so $d$ does not divide $ N_j.$

We conclude that $Z_{top, \I}^{(d)} = 0.$ \qed ${}$
\\ \\ \\
\emph{Acknowledgement:} The authors are grateful to Willem Veys for
the proposal to work on this problem.

\footnotesize{

}

\end{document}